\documentclass[12pt,a4paper]{amsart}
\usepackage{amsmath,amsfonts,amssymb,amscd,amsthm}
\theoremstyle{plain} %% This is the default
\newtheorem{Theorem}[equation]{Theorem}

\newtheorem{Lemma}[equation]{Lemma}
\newtheorem{Proposition}[equation]{Proposition}
\newtheorem{Fact}[equation]{Fact}

\theoremstyle{definition}

\newtheorem{Conjecture}[equation]{Conjecture}

\newtheorem{Assumption}[equation]{Assumption}

\theoremstyle{remark}
\newtheorem{Remark}[equation]{Remark}

\newtheorem{Example}[equation]{Example}

\numberwithin{equation}{section}
\numberwithin{figure}{section}
\newcommand{\R}{{\mathbb R}}
\newcommand{\C}{{\mathbb C}}

\newcommand{\posreal}{{\mathbb R}_{> 0}}

\newcommand{\Proj}[2]{{\mathbb P}^{#1}\!\left({#2}\right)}

\newcommand{\ric}[1]{\operatorname{Ric}_{#1}}

\newcommand{\cinf}{C^{\infty}}

\newcommand{\map}[3]{{#1}\colon{#2}\longrightarrow{#3}}
\begin{document}
\title{
New examples of Sasaki-Einstein manifolds
}
\author{Toshiki Mabuchi*}
\address{*Department of Mathematics, Osaka University,
Toyonaka, Osaka 560-0043, Japan}
\email{mabuchi@math.sci.osaka-u.ac.jp}
\author{Yasuhiro Nakagawa**}
\address{
**School of Mathematics and Physics, College of Science and Engineering,
Kanazawa University,
Kakuma-machi, Kanazawa, 920-1192, Japan
}
\email{yasunaka@kenroku.kanazawa-u.ac.jp}
\thanks{2010 {\it Mathematics Subject Classification.}
Primary~53C25;\, Secondary~32Q20, 53C55}
\thanks{
*Supported by JSPS Grant-in-Aid for Scientific Research (A) No.~20244005.
}
\thanks{
**Supported by JSPS Grant-in-Aid for Scientific Research
%(A) No.~20244005 and 
(C) No.~20540069.
}
%\subjclass{Primary 32Q20; Secondary 53C55, 53C25}
%
%\dedicatory{
%
%}
%%%%%%%%%%%%%%%%%%%%%%%%%%%%%%%%%%%%%%%%%%%%%%%%%%%%%%%%%%%%%%%%%%%%%%
\begin{abstract}
 In this note, stimulated by the existence results
 \cite{Futaki-Ono-Wang09a} %, \cite{Cho-Futaki-Ono08a}
 for toric Sasaki-Einstein metrics, we obtain 
new examples of Sasaki-Einstein metrics on 
 $S^1$-bundles associated to canonical line bundles of
 $\Proj{1}{\C}$-bundles over K\"{a}hler-Einstein Fano manifolds, even
 though the Futaki's obstruction does not vanish. 
 Here the method as in
 \cite{Sakane86a}, \cite{Koiso-Sakane86a}, \cite{Mabuchi87a} is used, and
our examples include non-toric Sasaki-Einstein manifolds.
\end{abstract}
\maketitle
\pagestyle{myheadings}
\markboth{TOSHIKI MABUCHI AND YASUHIRO NAKAGAWA}
{NEW EXAMPLES OF SASAKI-EINSTEIN MANIFOLDS}
%%%%%%%%%%%%%%%%%%%%%%%%%%%%%%%%%%%%%%%%%%%%%%%%%%%%%%%%%%%%%%%%%%%%%%
%% intro.tex (LaTeX2e file)
%%%%%%%%%%%%%%%%%%%%%%%%%%%%%%%%%%%%%%%%%%%%%%%%%%%%%%%%%%%%%%%%%%%%%%
\section{Introduction}
\label{section:intro}

Sasaki-Einstein manifolds were studied not only by mathematicians but
also by physicists, as Sasaki-Einstein manifolds have various
interesting phenomena such as ``AdS/CFT correspondence''
in theoretical physics
(cf.\ \cite{Boyer-Galicki99a}, \cite{Boyer-Galicki00a},
\cite{Boyer-Galicki08a}, \cite{Boyer-Galicki-Kollar05a},
\cite{Boyer-Galicki-Simanca08a},
%\cite{Futaki-Ono-Wang09a}, \cite{Cho-Futaki-Ono08a},
\cite{Hashimoto-Sakaguchi-Yasui04a}, \cite{Martelli-Sparks05a},
\cite{Martelli-Sparks06a}, \cite{Martelli-Sparks-Yau06a},
\cite{Martelli-Sparks-Yau08a}). Recently in  \cite{Cho-Futaki-Ono08a}
and \cite{Futaki-Ono-Wang09a}, classification of toric Sasaki-Einstein
manifolds was given.

\medskip
A Sasaki manifold is a $(2m+1)$-dimensional Riemannian manifold $(S,g)$
whose cone manifold $(C(S),\overline{g})$ is a K\"{a}hler manifold with
$$
C(S):= S\times\posreal
\;\;\text{ and }\;\;
 \overline{g}:=\left(dr\right)^2+r^2g,
 $$
where $r$ is the standard
coordinate on the set $\posreal =\{ r > 0\}$ of positive real
numbers. Then $S$ is a contact manifold with the contact form
$$
\eta:=\left(\sqrt{-1}\left(\overline{\partial}-\partial\right)\!
\log r\right)|_{r=1}.
$$
Here $S$ is viewed as the submanifold of
$C(S)$ defined by the equation $r=1$. We further consider the 
the {\it Reeb field} $\xi$ characterized by
$$
i(\xi)\eta=1 \,\text{ and }\,
i(\xi)d\eta=0,
$$
where $i(\xi)$ is the interior product by $\xi$.
The Reeb field $\xi$ is a Killing vector field on $(S,g)$ with a lift to
a holomorphic Killing vector field on $(C(S),\overline{g})$. This
generates a $1$-dimensional foliation on $S$,
called the {\it Reeb foliation}. The Sasaki metric $g$ naturally induces
a transverse K\"{a}hler metric $g^{\operatorname{T}}$ for the Reeb
foliation on $S$. A Sasaki manifold $(S,g)$ is {\it toric}, if $C(S)$ is
a toric manifold.

\medskip
The following well-known fact allows us to reduce the existence of
Sasaki-Einstein metrics to that of transverse K\"{a}hler-Einstein
metrics:

\begin{Fact}[cf.\ {\cite[Chapter~11]{Boyer-Galicki08a}}]
\label{fact:t-einstein}
 A Sasaki manifold $(S,g)$ is Einstein with positive scalar curvature
 $2m$ if and only if the transverse K\"{a}hler metric
 $g^{\operatorname{T}}$ is Einstein with positive scalar curvature
 $2(m+1)$.
\end{Fact}

We now pose the following conjecture:

\begin{Conjecture}\label{conj:conj}
 Let $M$ be a Fano manifold. If there exists a K\"{a}hler-Ricci soliton
 (see for instance \cite{Wang-Zhu04a} for K\"{a}hler-Ricci
 solitons) on $M$, then the $S^1$-bundle associated to the canonical
 line bundle $K_M$ of $M$ admits a Sasaki-Einstein metric with a
 suitable choice of the Reeb field.
\end{Conjecture}

By the results of Wang and Zhu \cite{Wang-Zhu04a}, the existence of
K\"{a}hler-Ricci solitons is known for toric Fano manifolds. Hence, the
results in \cite{Futaki-Ono-Wang09a} shows that
Conjecture~\ref{conj:conj} is affirmative for toric Fano manifolds.

\medskip
We now consider Koiso-Sakane's examples
\cite{Sakane86a}, \cite{Koiso-Sakane86a}, \cite{Koiso-Sakane88a} of
$\Proj{1}{\C}$-bundles 
over K\"{a}hler-Einstein Fano manifolds.
To fix our notation, recall the paper \cite{Mabuchi87a}.
Under the assumption below, we fix once for all a compact connected
$n$-dimensional complex manifold $W$ with $c_1(W)>0$ and an Hermitian
holomorphic line bundle $(L,h)$ over $W$.

\begin{Assumption}\label{assumption:tight-pair}
 (1)\,There exists a K\"{a}hler-Einstein form $\omega_0$ on $W$, i.e.,
 $\ric{}(\omega_0)=\omega_0$, where $\ric{}(\omega_0)$ is the Ricci form
 for $\omega_0$.

\smallskip\noindent
 (2)\,$2\pi c_1(L;h):=\sqrt{-1}\,\overline{\partial}{\partial}\log h$
 has constant eigenvalues
$$
\mu_1\leqq\mu_2\leqq\dotsm\leqq\mu_n
$$
 with respect to $\omega_0$ satisfying $-1<\mu_k<1$ for
 $k=1,2,\dotsc,n$.
\end{Assumption}\noindent
By this assumption, the compactification
$M_W^L:={\mathbb P}\left(L\oplus\mathcal{O}_W\right)$ of $L$ is 
%an $(n+1)$-dimensional complex manifold 
a $\Proj{1}{\C}$-bundle over $W$ with $c_1(M_W^L)>0$.
%The first named author proved the following:
%\begin{Fact}
%[Sakane~\cite{Sakane86a}; see also \cite{Mabuchi87a}]
%\label{fact:einstein-p1-bundle}
%Under Assumption~\ref{assumption:tight-pair} as above,
Then $M_W^L$ admits a K\"{a}hler-Einstein metric if and only if its
Futaki's obstruction (cf.\ \cite{Futaki83a}) vanishes:
\begin{equation}\label{eq:futaki-vanishing}
 \int_{-1}^{1}x\prod_{k=1}^n\left(1+\mu_kx\right)dx=0.
\end{equation}
%\end{Fact}

%\begin{Remark}
 %The condition \eqref{eq:futaki-vanishing} means nothing but the
 %vanishing of the Futaki character of $M_W^L$ (see  for
 %the definition of the Futaki character).
%\end{Remark}

\noindent
Let $S_W^L$ be the $S^1$-bundle over $M_W^L$ associated to the canonical
line bundle $K_{M_W^L}$ of $M_W^L$. In \cite{Koiso90a}, Koiso showed
that a K\"{a}hler-Ricci soliton exists on $M_W^L$, whether or not
equality~\eqref{eq:futaki-vanishing} holds. Hence by
Conjecture~\ref{conj:conj}, a Sasaki-Einstein metric is expected to
exist on $S_W^L$. The purpose of this note is to give the following
affirmative result:

\begin{Theorem}\label{thm:main1}
 Under the Assumption~\ref{assumption:tight-pair}, whether or not
 equality~\eqref{eq:futaki-vanishing} holds, $S_W^L$ always admits a
 Sasaki-Einstein metric for a suitable choice of the Reeb
 field. Furthermore, $K_{M_W^L}$ admits a complete Ricci-flat K\"{a}hler
 metric in every K\"{a}hler class.
\end{Theorem}

\begin{Remark}
 Kobayashi \cite{Kobayashi63a} (see also Jensen \cite{Jensen73a}, Wang
 and Ziller \cite{Wang-Ziller90a}) constructed Einstein metrics on
 $S^1$-bundles over Einstein manifolds. Our theorem
 %(Theorem~\ref{thm:main1})
 above shows that $S_W^L$ always admits an Einstein  metric, even though
 $M_W^L$ admits no K\"{a}hler-Einstein metrics.
\end{Remark}
%%%%%%%%%%%%%%%%%%%%%%%%%%%%%%%%%%%%%%%%%%%%%%%%%%%%%%%%%%%%%%%%%%%%%%
%% section2.tex (LaTeX2e file)
%%%%%%%%%%%%%%%%%%%%%%%%%%%%%%%%%%%%%%%%%%%%%%%%%%%%%%%%%%%%%%%%%%%%%%
\section{Transverse holomorphic structures on $S_W^L$}
\label{section:t-holo-str}

For an open cover $\left\{U_\alpha\,;\,\alpha\in A\right\}$ of $W$,
we choose  a system
 of holomorphic local
coordinates $(w_\alpha^1,w_\alpha^2,\dotsc,w_\alpha^n)$ on each
$U_\alpha$, and by taking a holomorphic local frame $e_\alpha$ for $L$,
we have the fiber coordinate $\zeta_\alpha^{+}$ for $L$ over $U_\alpha$
with respect to $e_\alpha$. Then
$(w_\alpha^1,w_\alpha^2,\dotsc,w_\alpha^n;\zeta_\alpha^{+})$ forms a
system of holomorphic local coordinates for
$U_\alpha^{+}:=L|_{U_\alpha}$. Let
$f_\alpha$ be the  frame for $L^{-1}$ dual to $e_\alpha$, and
let $\zeta_\alpha^{-}$ be the fiber coordinate for $L^{-1}$ over
$U_\alpha$ with respect to $f_\alpha$. Then
$(w_\alpha^1,w_\alpha^2,\dotsc,w_\alpha^n;\zeta_\alpha^{-})$ form a
system of holomorphic local coordinates on
$U_\alpha^{-}:=L^{-1}|_{U_\alpha}$. Then $U^{+}_{\alpha}$ and
$U^{-}_{\alpha}$ are glued together by the relation
$$
\zeta_\alpha^{+}=\left(\zeta_\alpha^{-}\right)^{-1}
$$
to form
$M_W^L=\mathbb{P}(L\oplus\mathcal{O}_W)=
\bigcup_{\alpha\in A}(U_\alpha^+\cup U_\alpha^-)$. Here,
$$
\pm dw_\alpha^1\wedge dw_\alpha^2\wedge\dotsm\wedge
dw_\alpha^n\wedge d\zeta_\alpha^\pm
$$
is a holomorphic local frame for $K_{M_W^L}$ over $U_\alpha^\pm$,
and with respect to this local frame, we have the fiber coordinate
$\tau_\alpha^\pm$ for $K_{M_W^L}$, respectively, i.e., all $(+)$-signs
and all $(-)$-signs should be chosen respectively.
Note that
{\allowdisplaybreaks\begin{align*}
 &\tau_\alpha^+
 dw_\alpha^1\wedge dw_\alpha^2\wedge\dotsm\wedge dw_\alpha^n\wedge
 d\zeta_\alpha^+\\&=
 \tau_\beta^+
 dw_\beta^1\wedge dw_\beta^2\wedge\dotsm\wedge dw_\beta^n\wedge
 d\zeta_\beta^+\\&=
 \tau_\beta^+ \phi_{\beta\alpha}(w)\psi_{\beta\alpha}(w)^{-1}
 dw_\alpha^1\wedge dw_\alpha^2\wedge\dotsm\wedge dw_\alpha^n\wedge
 d\zeta_\alpha^+,
\end{align*}}\noindent
for $w\in U_\alpha\cap U_\beta$.
Here $\{\psi_{\beta\alpha}\,;\,\alpha,\beta\in A\}$
are the transition functions for $L$ with respect to the local
frames $\{e_\alpha\,;\,\alpha\in A\}$ for $L$, 
while $\{\phi_{\beta\alpha}\,;\,\alpha,\beta\in A\}$ 
are the transition functions for $K_W$ with respect to the local
frames $\{dw_\alpha^1\wedge
 \dotsm\wedge dw_\alpha^n\,;\,\alpha\in A\}$ for $K_W$, i.e.,
\begin{align*}
 &e_\beta=\psi_{\beta\alpha}(w)e_\alpha,\quad
 \quad f_\beta=\psi_{\beta\alpha}(w)^{-1} f_\alpha,\\
 &dw_\beta^1\wedge dw_\beta^2\wedge \dotsm\wedge dw_\beta^n
  =\phi_{\beta\alpha}(w)
  dw_\alpha^1\wedge dw_\alpha^2\wedge \dotsm\wedge dw_\alpha^n,
\end{align*}
for $w\in U_\alpha\cap U_\beta$.
Hence $\tau_\alpha^+$ can be viewed as the fiber coordinate for
$K_W\otimes L^{-1}$ over $U_\alpha$ with respect to the local frame
$\left(dw_\alpha^1\wedge\dotsm\wedge
 dw_\alpha^n\right)\otimes f_\alpha$. Similarly, $\tau_\alpha^-$ is also
viewed as the fiber coordinate for $K_W\otimes L$ over $U_\alpha$ with
respect to the local frame
$\left(dw_\alpha^1\wedge\dotsm\wedge dw_\alpha^n\right) \otimes
e_\alpha$. Moreover, since
$\tau_{\alpha}^+\zeta_\alpha^+ =\tau_{\alpha}^-\zeta_\alpha^-$
on $U_\alpha^+\cap U_\alpha^-$, it follows that 
$$
\tau_\alpha^+(\zeta_\alpha^+)^2=\tau_\alpha^-.
$$
Now, for $-\frac{1}{2}< a\in \mathbb R$, we consider holomorphic vector
fields
{\allowdisplaybreaks\begin{align*}
&a\sqrt{-1}\zeta_\alpha^+\frac{\partial}{\partial\zeta_\alpha^+}
  +\sqrt{-1}\tau_\alpha^+\frac{\partial}{\partial\tau_\alpha^+}
  \,\,\,\,\text{ on }\tilde{p}^{-1}(U_\alpha^+),\\
&\!-a\sqrt{-1}\zeta_\alpha^-\frac{\partial}{\partial\zeta_\alpha^-}
  +(1+2a)\sqrt{-1}\tau_\alpha^-\frac{\partial}{\partial\tau_\alpha^-}
  \,\,\,\,\text{ on }\tilde{p}^{-1}(U_\alpha^-),
\end{align*}}\noindent
where $\map{\tilde{p}}{K_{M_W^L}}{M_W^L}$ is the natural projection.
Then these are glued together to define a well-defined global
holomorphic vector field $\xi_a$ on $K_{M_W^L}$. We choose
$\xi_a+\bar{\xi}_a$ as the Reeb field on $S_W^L$. However, we
call $\xi_a$ also as the Reeb filed by abuse of terminology. Put
$$
z_\alpha^+:=(\tau_\alpha^+)^{-a}\zeta_\alpha^+ \;\; \text{ and }\;\;
z_\alpha^-:=(\tau_\alpha^-)^{\frac{a}{1+2a}}\zeta_\alpha^-. 
$$
Then
$(w_\alpha^1,w_\alpha^2,\dotsc,w_\alpha^n;z_\alpha^+)$ and
$(w_\alpha^1,w_\alpha^2,\dotsc,w_\alpha^n;z_\alpha^-)$ are transverse
holomorphic local coordinates on
$\widetilde{U}_\alpha^+:=p^{-1}(U_\alpha^+)$ and
$\widetilde{U}_\alpha^-:=p^{-1}(U_\alpha^-)$, respectively,
with respect to the Reeb field $\xi_a$, in view of the identities
$$
dz_\alpha^+(\xi_a)=0\;\;\text{ and }\;\;
dz_\alpha^-(\xi_a)=0,
$$
where $\map{p}{S_W^L}{M_W^L}$ is the natural projection.
Note that $z_\alpha^+$ and $z_\alpha^-$ satisfy the
following relation:
$$
z_\alpha^+=(\tau_\alpha^+)^{-a}\zeta_\alpha^+=
(\tau_\alpha^-)^{-a}(\zeta_\alpha^-)^{-(1+2a)}=
(z_\alpha^-)^{-(1+2a)}.
$$
For the natural projection $\map{q}{S_W^L}{W}$, the fiber $q^{-1}(w)$ 
over each $w\in U_\alpha$ has a transverse holomorphic structure defined
by the transverse holomorphic coordinate $z_\alpha^\pm$. Then on
$q^{-1}(w)$,
$$
G:=
\begin{cases}
\;\; \left(|z_\alpha^+|^{-1}+|z_\alpha^+|^{\frac{1}{1+2a}}\right)^{\!\!-2}
 \dfrac{|dz_\alpha^+|^2}
       {|z_\alpha^+|^{2}},&
 \text{ on }\,\,q^{-1}(w)\cap\widetilde{U}_\alpha^+,\\
\;\; \left(1+2a\right)^{2}
 \left(|z_\alpha^-|^{1+2a}+|z_\alpha^-|^{-1}\right)^{-2}
 \dfrac{|dz_\alpha^-|^2}
       {|z_\alpha^-|^{2}},&
 \text{ on }\,\,q^{-1}(w)\cap\widetilde{U}_\alpha^+,
\end{cases}
$$
defines a transverse K\"{a}hler metric, which is invariant under the
standard $S^1$-action
$z_\alpha^+\overset{t}{\longmapsto}tz_\alpha^+$, $t\in S^1$, for
each $w\in U_\alpha$. By setting $x:=-2\log|z_\alpha^+|$, we define
$$
v(x):=2\log\left\{ \exp \left (\frac{x}{2}\right )+\exp \left
({-\frac{x}{2(1+2a)}}\right )\right\}.
$$
Then its derivative $v^\prime(x)$ defines a moment map whose image is
the closed interval $\left[-\frac{1}{1+2a},\,1\right]$.

%%%%%%%%%%%%%%%%%%%%%%%%%%%%%%%%%%%%%%%%%%%%%%%%%%%%%%%%%%%%%%%%%%%%%%
%% section3.tex (LaTeX2e file)
%%%%%%%%%%%%%%%%%%%%%%%%%%%%%%%%%%%%%%%%%%%%%%%%%%%%%%%%%%%%%%%%%%%%%%
\section{Sasaki-Einstein metrics on $S_W^L$}
\label{section:t-ke}

In this section, we construct an Sasaki-Einstein metric on $S_W^L$
by an argument as in \cite{Mabuchi87a}.
For $a>-\frac{1}{2}$, define a polynomial $A_a(x)$ in $x$ by
$$
A_a(x):=-\int_{-\frac{1}{1+2a}}^{x}s\prod_{k=1}^n
\left (1+\mu^{}_{k,a} s\right )ds.
$$
where $\mu_{k,a}:=\mu_k+a(1+\mu_k)$ for $k=1,2,\dotsc,n$.
Now, we assume that $A_a(1)=0$. Since $a>-\frac{1}{2}$, it follows from
Assumption~\ref{assumption:tight-pair} that
{\allowdisplaybreaks\begin{align*}
 &0<A_a(x)\leqq A_a(0),\\
 &\frac{A_a^\prime(x)}{x}<0,
\end{align*}}\noindent
for $-\frac{1}{1+2a}<x<1$. In particular, the rational function
$\frac{A_a^\prime(x)}{xA_a(x)}$ is free from poles and zeros over the
open interval $\left(-\frac{1}{1+2a},\,1\right)$ and has a pole of order
$1$ at both $x=-\frac{1}{1+2a}$ and $x=1$. Hence,
$$
B_a(x):=-\int_{0}^{x}\frac{A_a^\prime(s)}{sA_a(s)}ds
$$
is monotone increasing over the interval
$\left(-\frac{1}{1+2a},\,1\right)$ and moreover, $B_a$ maps
$\left(-\frac{1}{1+2a},\,1\right)$ diffeomorphically onto $\R$. Let
$$
\map{B_a^{-1}}{\R}{\left(-\frac{1}{1+2a},\,1\right)}
$$ 
be the inverse function of
$\map{B_a}{\left(-\frac{1}{1+2a},\,1\right)}{\R}$, and define a
$\cinf$ functions $x_a(\rho)$ and $u_a(\rho)$ in $\rho\in\R$ by
$x_a(\rho):=B_a^{-1}(\rho)$ and $u_a(\rho):=-\log(A_a(x_a(\rho)))$.
Then $u_a^\prime(\rho)=x_a(\rho)$ and hence
\begin{equation}\label{eq:diff-eq}
 u_a^{\prime\prime}(\rho)\prod_{k=1}^{n}
 \left (1+\mu^{}_{k,a}u_a^\prime(\rho)\right )
 =e^{-u_a(\rho)}.
\end{equation}
%where $\mu_{k,a} := \mu_k + a(1+\mu_k )$. Now, 
On $\widetilde{U}_\alpha^+$, we define
$\rho_\alpha^+:=-\log|z_\alpha^+|^2
-\log\left(\kappa_\alpha^{-a}h_\alpha^{1+a}\right)$ by setting
{\allowdisplaybreaks\begin{align*}
 &\kappa_\alpha:=h_{K_W}\!\!
  \left(dw_\alpha^1\wedge\dotsm\wedge
        dw_\alpha^n,
        dw_\alpha^1\wedge\dotsm\wedge
        dw_\alpha^n\right),\\
 &h_\alpha:=h(e_\alpha,e_\alpha),
\end{align*}}\noindent
that is, $\exp \left ({-\frac{\rho_\alpha^+}{2}}\right )$ can be
formally viewed as the norm of
$$
z_\alpha^+\left(
 \left(\frac{\partial}{\partial w_\alpha^1}\wedge\dotsm\wedge
       \frac{\partial}{\partial w_\alpha^n}\right)^{\!\!a}
       \otimes e_\alpha^{1+a}\right)
$$
with respect to the Hermitian metric
$h_{K_W}^{-a}\otimes h^{1+a}$ on $K_W^{-a}\otimes L^{1+a}$.
Here $h_{K_W}$ denotes the Hermitian metric for $K_W$ induced
by $\omega_0$. Then we have $\rho_\alpha^+=\rho_\beta^+$ on
$\widetilde{U}_\alpha^+\cap\widetilde{U}_\beta^+$. Now we consider the
following transverse $(n+1,n+1)$-form $\varPhi_\alpha$, with respect to
$\xi_a$, on $\widetilde{U}_\alpha^+$:
$$
\varPhi_{\alpha,+}:=\sqrt{-1}\,(n+1)\exp ({-u_a(\rho_\alpha^+)})
\left(q^*\omega_0\right)^n\wedge
\frac{dz_\alpha^+\wedge d\overline{z_\alpha^+}}{|z_\alpha^+|^2}.
$$
Note that
$\ric{}(\omega_0)=\sqrt{-1}\,\overline{\partial}\partial\log\omega_0^n$
and that, for each fixed $w_0\in U_\alpha$, we can choose a local frame
$e_\alpha$ for $L$ and a system
$(w_\alpha^1,w_\alpha^2,\cdots,w_\alpha^n)$ of holomorphic local
coordinates on $U_\alpha$ satisfying
{\allowdisplaybreaks\begin{align*}
 &d\left(\kappa_\alpha^{-a}h_\alpha^{1+a}\right)(w_0)=0,\\
 &\omega_0(w_0)=
  \sqrt{-1}\sum_{k=1}^{n}dw_\alpha^k\wedge d\overline{w_\alpha^k},\\ 
 &\!\left(\sqrt{-1}\,\overline{\partial}\partial\log h_\alpha\right)
  (w_0)=\sqrt{-1}\sum_{k=1}^{n}
  \mu_k dw_\alpha^k\wedge d\overline{w_\alpha^k}.
\end{align*}}\noindent
Then, along $q^{-1}(w_0)\cap\widetilde{U}_\alpha^+$, 
%by setting  $\mu_{k,a}:= \mu_k + a (1+\mu_k )$, 
we write
 ${\omega}^{\operatorname{T}}_{\alpha,+}:=
 \sqrt{-1}\,\overline{\partial}\partial\log\varPhi_{\alpha,+}$ 
as a sum
 $$
 %{\omega}^{\operatorname{T}}_{\alpha,+} = 
 %{\allowdisplaybreaks\begin{align*}
% {\omega}^{\operatorname{T}}_{\alpha,+}:=&
% \sqrt{-1}\,\overline{\partial}\partial\log\varPhi_{\alpha,+}\\=&
 u^{\prime\prime}(\rho_\alpha^+)
 \frac{\sqrt{-1}dz_\alpha^+\wedge d\overline{z_\alpha^+}}
      {|z_\alpha^+|^2}
%\\
%&
 +\sum_{k=1}^{n}\left\{
 \left (1+\mu_{k,a} u^\prime(\rho_\alpha^+)\right )
 \sqrt{-1}dw_\alpha^k\wedge d\overline{w_\alpha^k}\right\}.
 $$
%\end{align*}}
\noindent
Since $a>-\frac{1}{2}$ and $-1<\mu_k<1$ ($k=1,2,\dotsc,n$),
${\omega}^{\operatorname{T}}_\alpha$ is a transverse K\"{a}hler form,
with respect to $\xi_a$, on
$\widetilde{U}_\alpha^+\setminus\{z_\alpha^+=0\}$. Furthermore, 
by \eqref{eq:diff-eq}, we have
$\left({\omega}^{\operatorname{T}}_{\alpha,+}\right)^{n+1}
 =\varPhi_{\alpha,+}$.
Therefore, ${\omega}^{\operatorname{T}}_\alpha$ defines a transverse
K\"{a}hler-Einstein metric, with respect to $\xi_a$, on
$\widetilde{U}_\alpha^+\setminus\{z_\alpha^+=0\}$. By an argument
as in \cite{Mabuchi87a}, the condition $A_a(1)=0$ implies
that
$\left\{{\omega}^{\operatorname{T}}_{\alpha,+}\,;,\alpha\in A\right\}$
are glued together to define a well-defined global transverse
K\"{a}hler-Einstein form ${\omega}^{\operatorname{T}}$ on $S_W^L$ with
the Reeb field $\xi_a$. 

\begin{Remark}
 $\varPhi_{\alpha,+}$ is formally viewed as an Hermitian metric on
 $K_W^{-1}\otimes\left(K_W^a\otimes L^{-(1+a)}\right)^{-1}=
 \left(K_W\otimes L^{-1}\right)^{-(1+a)}$.
 Then we put
{\allowdisplaybreaks\begin{align*}
 &r:=|\tau_\alpha^+|\left \{ (n+1)\,
  \exp(-u(\rho_\alpha^+))\,\kappa_\alpha^{-1}|z_\alpha^+|^{-2} 
  \right \}^{\!\!-\frac{1}{2(1+a)}},\\
 &\eta:=
  \left(\sqrt{-1}(\overline{\partial}-\partial)\log r\right)|_{r=1}.
\end{align*}}\noindent
 Since $r = r (\tau_\alpha^+)$ is regarded as the norm of
$$
\tau_\alpha^+\left( 
\left(\frac{\partial}{\partial w_\alpha^1}\wedge\dotsm\wedge
      \frac{\partial}{\partial w_\alpha^n}\right)\otimes
\frac{\partial}{\partial z_\alpha^+}\right)^{\!\!-\frac{1}{1+a}},
$$
 with respect to the Hermitian metric
 $\left(\varPhi_{\alpha,+}\right)^{-\frac{1}{1+a}}$ for
 $K_W\otimes L^{-1}$, $r$ defines a well-defined
 $\cinf$ function on $K_{M_W^L}$, and in particular $S_W^L$ 
 is identified with the submanifold of $K_{M_W^L}$ defined by the
 equation $r=1$. Moreover,
$$
g:=(\eta)^2+\frac{1}{1+a}g^{\operatorname{T}}
$$ 
 is a Riemann metric on $S_W^L$ and $\eta$ is a contact form on $S_W^L$,
 where $g^{\operatorname{T}}$ is the transverse K\"{a}hler metric
 associated to $\omega^{\operatorname{T}}$. Furthermore, the fundamental
 form $\overline{\omega}$ of the cone metric $\overline{g}$ associated
 to $g$ is given by
$$
\overline{\omega}:=
2rdr\wedge\eta+\frac{r^2}{1+a}\omega^{\operatorname{T}}.
$$
 In view of $d\eta=\frac{1}{1+a}\omega^{\operatorname{T}}$, we obtain
 $d\overline{\omega}=0$, and hence $(S_W^L,g)$ is a Sasaki manifold with
 the Reeb field $\xi_a$. 
\end{Remark}

Now by Fact~\ref{fact:t-einstein}, we obtain the following criterion on
the existence of Sasaki-Einstein metrics on $S^L_W$:
\begin{Proposition}\label{prop:prop1}
 Under the Assumption~\ref{assumption:tight-pair}, if
\begin{equation}\label{eq:soliton-field}
 A_a(1)=-\int_{-\frac{1}{1+2a}}^1x
 \prod_{k=1}^n\left (1+\mu_{k,a}x \right ) dx=0,
\end{equation}
 then $S_W^L$ admits a Sasaki-Einstein metric with the Reeb field
 $\xi_a$.
\end{Proposition}
\begin{Remark}
 In the special case $a=0$, we easily see that \eqref{eq:soliton-field}
 is nothing but the condition~\eqref{eq:futaki-vanishing} in the
 introduction.
\end{Remark}
Next, we shall show the existence of $a\in \R$  such that both
$a>-\frac{1}{2}$ and $A_a(1)=0$ hold. We now put
{\allowdisplaybreaks\begin{align*}
 &f(x;a):=x\prod_{k=1}^{n}
  \left (1+\mu_{k,a} x\right ),\\
 &F(a):=\int_{-\frac{1}{1+2a}}^1f(x;a)\, dx\,\,\,\,(=-A_a(1)).
\end{align*}}\noindent
%Since $F_W^{L^{-1}}(0)=-F_W^L(0)$ and $S_W^{L^{-1}}=S_W^L$, we may
%assume that $F_W^L(0)\leqq 0$ without loss of generality. Since
%$\lim_{a\to\infty}\int_0^1f(x;a)dx=+\infty$, we have
%$$
%\lim_{a\to\infty}F_W^L(a)=+\infty.
%$$
%In view of the continuity of $F_W^L$, we can find $a_0>0$ such that
%$F_W^L(a_0)=0$.
Since $\lim_{a\to+\infty}F(a)=+\infty$ and
$\lim_{a\to-\frac{1}{2}+0}F(a)=-\infty$,
the continuity of $F$ allows us to find
$a_0>-\frac{1}{2}$ such that $F(a_0)=0$.
Moreover, 
$$
{F}^\prime(a)=
\int_{-\frac{1}{1+2a}}^1\frac{\partial}{\partial a}f(x;a)dx +
\frac{-2}{(1+2a)^2}f\left(-\frac{1}{1+2a};a\right).
$$
Note also that $\mu_{k,a} = \mu_k + a (1+\mu_k )$. Hence for
$-\frac{1}{1+2a}\leqq x\leqq 1$,
{\allowdisplaybreaks\begin{align*}
 &\frac{\partial}{\partial a}f(x;a)=x^2\sum_{j=1}^{n}\left\{(1+\mu_j)
  \prod_{k\ne j}
  \left (1+\left\{\mu_k+a(1+\mu_k)\right\}x\right)\right\}
  \geqq 0,\\
 &f\left(-\frac{1}{1+2a};a\right)=
  -\left(\frac{1}{1+2a}\right)^{\!\! n+1}
  \prod_{k=1}^{n}\left\{(1+a)(1-\mu_k)\right\}<0.
\end{align*}}\noindent
Now in the expression of $F'(a)$, the first term is nonnegative and the
second term is positive. Therefore ${F}^\prime(a)>0$. Hence we obtain:
\begin{Lemma}\label{lem:reeb-field}
 Under the Assumption~\ref{assumption:tight-pair}, there exists a unique
 real number $a_0>-\frac{1}{2}$ such that $F(a_0)=0$.
\end{Lemma}

Therefore, by Proposition~\ref{prop:prop1} and
Lemma~\ref{lem:reeb-field}, if Assumption~\ref{assumption:tight-pair}
is satisfied, then $S_W^L$ always admits an Sasaki-Einstein metric with
the Reeb field $\xi_{a_0}$. On the other hand, in view of
\cite{Futaki10a}, \cite{Goto09a}(see also \cite{vanCoevering08a}),
we now conclude that $K_{M_W^L}$ admits a complete Ricci-flat K\"{a}hler
metric in every K\"{a}hler class. The proof of
Theorem~\ref{thm:main1} is now complete.

%%%%%%%%%%%%%%%%%%%%%%%%%%%%%%%%%%%%%%%%%%%%%%%%%%%%%%%%%%%%%%%%%%%%%%
%% section4.tex (LaTeX2e file)
%%%%%%%%%%%%%%%%%%%%%%%%%%%%%%%%%%%%%%%%%%%%%%%%%%%%%%%%%%%%%%%%%%%%%%
\section{Examples}
\label{section:examples}

In this section, we shall give a couple of examples of Sasaki-Einstein
manifolds as an application of Theorem~\ref{thm:main1}.

\begin{Example}\label{example:ex1}
 We first put
\begin{align*}
 &W:=\prod_{i=1}^{l}\Proj{\,n_i}{\C},\\
 &L:=\bigotimes_{i=1}^{l}
  p_i^*\left(\mathcal{O}_{\Proj{n_i}{\C}}(\nu_i)\right),
\end{align*}
 where $\map{p_i}{W}{\Proj{\,n_i}{\C}}$ is the natural projection to
 the $i$-th factor ($i=1,2,\dotsc,l$).
 In view of the isomorphism
 $K_{\Proj{k}{\C}}^{-1}\cong\mathcal{O}_{\Proj{k}{\C}}(k+1)$, if
%$n_i+1\pm\nu_i>0,\,\,(i=1,2,\dotsc,l)$,
$$
-(n_i+1)<\nu_i<n_i+1,\qquad(i=1,2,\dotsc,l),
$$
 then the pair $(W,L)$ satisfies
 Assumption~\ref{assumption:tight-pair}. Hence by
 Theorem~\ref{thm:main1}, $S_W^L$ admits a Sasaki-Einstein metric,
 though this is toric. Then $F(a)$ in Section 3 is given by
$$
F(a)=\int_{-\frac{1}{1+2a}}^{1}
x\prod_{i=1}^{l}\left ( 1+\left\{\frac{\nu_i}{n_i+1}+
a\left(1+\frac{\nu_i}{n_i+1}\right)\right\}x\right )^{n_i}dx.
$$
 For instance, we consider the simplest case, that is, $W=\Proj{1}{\C}$
 and $L=\mathcal{O}_{\Proj{1}{\C}}(1)$. In this case, $M_W^L$ is a del
 Pezzo surface obtained from $\Proj{2}{\C}$ by blowing up one point, and
 we see the irregularity of $(S_W^L,\xi_{a_0})$ by
$$
a_0=\frac{-5+\sqrt{13}}{12}.
$$
\end{Example}

\begin{Example}\label{example:ex2}
 Next, let $W:=\operatorname{Gr}(k,p)$ be the complex Grassmannian
 manifold of all $p$-dimensional subspaces of $\C^k$, which is a complex
 manifold of dimension $p(k-p)$. Then there exists an ample
 line bundle $A(k,p)$ over $\operatorname{Gr}(k,p)$ such that
 $K_{\operatorname{Gr}(k,p)}^{-1}\cong A(k,p)^{k}$ (see for instance
 \cite[p.\ 205]{Takeuchi78a}). We put  $L:=A(k,p)^\nu$. If %$n\pm\nu>0$,
 $-k<\nu<k$, then the pair $(W,L)$ satisfies
 Assumption~\ref{assumption:tight-pair}. Hence by
 Theorem~\ref{thm:main1}, $S_W^L$ admits a Sasaki-Einstein metric, and
 if $2 \leqq p\leqq k-2$, then $S_W^L$ is non-toric.
\end{Example}

\begin{Example}
 Let $\mathcal{M}_n$ be the moduli space of smooth hypersurfaces of
 degree $n$ in $\Proj{n+1}{\C}$. For the Fermat type hypersurface
$$
W_0:=\left\{
[X_0:\dotsc:X_{n+1}]\,;\,
\left(X_0\right)^n+\dotsm+\left(X_{n+1}\right)^n=0
\right\}\in\mathcal{M}_n,
$$
 a theorem of Tian \cite{Tian87a} shows that $W_0$ admits a
 K\"{a}hler-Einstein metric, and in particular
$$
\mathcal{M}_n^{\operatorname{KE}}:=\left\{W\in\mathcal{M}_n\,;\,
W\text{ admits a K\"{a}hler-Einstein metric.}\right\}
$$
 is a non-empty open subset of $\mathcal{M}_n$. For every
$W\in\mathcal{M}_n^{\operatorname{KE}}$,
 we have $K_W\cong\mathcal{O}_{\Proj{n+1}{\C}}(-2)|_W$ by adjunction
 formula. Put $L:=\mathcal{O}_{\Proj{n+1}{\C}}(1)|_W$. Then the pair
 $(W,L)$ satisfies Assumption~\ref{assumption:tight-pair}, and
 Theorem~\ref{thm:main1} shows that $S_W^L$ admits a Sasaki-Einstein
 metric. If $n=3$, $W$ is a well-known cubic threefold, and in this case
 by \cite[Theorem~13.12]{Clemens-Griffiths72a}, $W$ is not birationally
 equivalent to $\Proj{3}{\C}$, and $S^L_W$ is again non-toric.
\end{Example}

%%%%%%%%%%%%%%%%%%%%%%%%%%%%%%%%%%%%%%%%%%%%%%%%%%%%%%%%%%%%%%%%%%%%%%
%% ref.tex (LaTeX2e file)
%%%%%%%%%%%%%%%%%%%%%%%%%%%%%%%%%%%%%%%%%%%%%%%%%%%%%%%%%%%%%%%%%%%%%%

%%%%%%%%%%%%%%%%%%%%%%%%%%%%%%%%%%%%%%%%%%%%%%%%%%%%%%%%%%%%%%%%%%%%%%
\end{document}